\documentclass[11pt]{amsart}
\usepackage{pxfonts}
\usepackage{enumitem}
\usepackage[T1]{fontenc}

\usepackage[pdftitle= {Kawamata },pdfauthor={MMLRPST}]{hyperref}
\usepackage{amscd, amssymb, amsmath, wasysym, yfonts}
\usepackage{graphicx}
\usepackage[all]{xy}
\usepackage{url}
\usepackage[usenames, dvipsnames]{color}
\usepackage{tcolorbox}
\usepackage{hyperref}
 \hypersetup{
  colorlinks=true,
   citecolor=purple!60!black!49!blue,
   linkcolor=black,
   urlcolor=teal}
\makeatletter
\def\l@subsection{\@tocline{2}{0pt}{2.5pc}{5pc}{}}
\makeatother
   
   \usepackage{enumitem}

\usepackage{amsrefs}
\usepackage[margin=1.50in]{geometry}
\theoremstyle{plain}
\newtheorem{theorem}{Theorem}[section]
\newtheorem*{thm*}{Theorem}
\newtheorem{thm}{Theorem}

\newtheorem{introcor}[thm]{Corollary}
\newtheorem*{thmB}{Theorem A*}
\newtheorem{prop}[theorem]{Proposition}

\newtheorem{lm}[theorem]{Lemma}

\theoremstyle{definition}

\newtheorem*{ex*}{Example}

\newcommand\OO{{\mathcal O}}

\newcommand\ol{\overline}

\DeclareMathOperator{\Pic}{Pic}
\DeclareMathOperator{\Proj}{Proj}
\DeclareMathOperator{\Aut}{Aut}
\DeclareMathOperator{\Diag}{Diag}

\newcommand\pp{{\mathbb{P}}}

\usepackage[colorinlistoftodos]{todonotes}

\title
{A footnote to a
theorem of Kawamata}
\author{Margarida Mendes Lopes}
\address{CAMGSD/Departamento de Matem\'atica,   Instituto Superior T\'ecnico, Universidade de Lisboa,  Av. Rovisco Pais, 1049-001 Lisboa, Portugal  }
\email{\url{mmendeslopes@tecnico.ulisboa.pt}}
\author{Rita Pardini}
\address{Dipartimento di Matematica, Universit\`{a} degli studi di Pisa, Largo Pontecorvo 5, 56127 Pisa (PI), Italy}
\email{\url{rita.pardini@unipi.it}}

\author{Sofia Tirabassi}

\address{Department of Mathematics, Stockholm University, Albano Campus, Stockholm, Sweden}
\email{\url{tirabassi@math.su.se}}

\newlist{todolist}{itemize}{2}
\setlist[todolist]{label=$\square$}
\usepackage{pifont}
\newcommand{\cmark}{\ding{51}}%
\newcommand{\done}{\rlap{$\square$}{\raisebox{2pt}{\large\hspace{1pt}\cmark}}%
\hspace{-2.5pt}}

\newlist{donelist}{itemize}{2}
\setlist[donelist]{label=\done}
\begin{document}

\begin{abstract}

   Kawamata has shown that the quasi-Albanese map of a quasi-projective  variety  with log-irregularity equal to the dimension  and log-Kodaira dimension 0 is birational.  In this note we show that under these hypotheses the quasi-Albanese map is proper in codimension 1 as conjectured by Iitaka.   

\end{abstract}

 \maketitle
 
\section*{Introduction}

Given    a smooth complex quasi-projective  variety  $V$, by Hironaka's theorem on the resolution of  singularities we can embedd $V$ into a smooth projective variety $X$ such that the complement $D:=X\backslash V$ is a reduced divisor on $X$ with simple normal crossings.  \par 
 One can use the compactification $X$ to define the logarithmic invariants of $V$. In particular, if $K_X$ denotes the canonical divisor  of $X$, then:
\begin{itemize}[label=-]

 \item the \emph{log-Kodaira dimension} of $V$ is $\overline{\kappa}(V):=\kappa(X,K_X+D);$
 \item the \emph{log-irregularity} of $V$ is $\overline{q}(V):=h^0(X,\Omega^1_X(\log D))$.
\end{itemize}

It easy to show that these invariants do not depend on the choice of the compactification $X$.\par
In a similar way to   projective varieties,  we can associate  to $V$ a quasi-abelian variety (i.e., an algebraic group which does not contain $\mathbb{G}_a$) $A(V)$, called the \emph{ quasi-Albanese variety} of $V$ and a  morphism $a_V\colon V\rightarrow A(V)$ which is called the \emph{quasi-Albanese morphism}. 

Kawamata, in his celebrated work \cite{Ka81}, provides a numerical criterion for the birationality of the quasi-Albanese map:
\begin{thm*}[Kawamata, Theorem 28 of \cite{Ka81}]
Let $V$ be a smooth complex algebraic variety of dimension $n$. If \ $\overline{\kappa}(V)=0$ and   $\overline{q}(V)=n$ then the quasi-Albanese morphism $a_V\colon V\rightarrow A(V)$ is birational.
\end{thm*}
In this note, we strengthen Kawamata's Theorem in a way that was conjectured by Iitaka in \cite[p. 501]{Ii77}. Our main result is the following:

\begin{thm}\label{thm:main}
 Let $V$ be a smooth complex algebraic variety of dimension $n
 $ with  $\overline{\kappa}(V)=0$ and   $\overline{q}(V)=n$. Denote by $a_V\colon V\rightarrow A(V)$ the quasi-Albanese morphism (which is birational).
 
  Then there is a closed subset $Z\subset A(V)$ of codimension $>1$ such that, setting $V^0:=V\backslash a_V^{-1}(Z)$,  the restriction
$a_{V|_{V^0}}\colon {V^0}\rightarrow A(V)\backslash Z$ is proper. 
\end{thm}

In  \cite {Ii77} Iitaka defined the quite involved notion of   WWPB (``weakly weak proper birational'')  equivalence. With  this language,  we can rewrite our statement in the following  way:
\begin{thmB}
 Let $V$ be a smooth complex algebraic variety of dimension $n
 $. Then $V$ is WWPB equivalent to a quasi-abelian variety if  and only if   $\overline{\kappa}(V)=0$ and   $\overline{q}(V)=n$. 
\end{thmB}
While the usual Kodaira dimension and irregularity are birational invariants, this is not the case for logarithmic invariants. As an easy example of this phenomenon, one can think of $\mathbb{P}^1$ and $\mathbb{P}^1\backslash\{p_0+p_1\}$: they are birationally equivalent but their invariants are completely different. However, WWPB-equivalent varieties share the same logarithmic invariants. Therefore this seems  to be the right notion of equivalence to consider when studying the birational geometry of open varieties. 

We conclude this section by remarking that the condition of being WWPB-equivalent is a really strong one: WWPB maps between normal affine varieties are indeed isomorphisms. Thus we get the following corollary of Theorem A:
\begin{introcor}
 A smooth complex affine variety $V$ of dimension $n$ is isomorphic to $\mathbb{G}_m^n$ if  and only if  $\overline{\kappa}(V)=0$ and $\overline{q}(V)=n$.
\end{introcor}

\noindent{\bf Notice.} Our original proof contained an unsubstantiated claim that was pointed to us by Prof. Osamu Fujino. In the attached addendum/correction (in collaboration with Prof. Fujino) we provide a new proof and we also explain how to avoid the gap in the original argument.

\noindent{\bf Notation.} 
We work over the complex numbers. If $X$ is  a smooth projective  variety, we denote by  $K_X$ the canonical class, and by $q(X):=h^0(X,\Omega^1_X)=h^1(X,\mathcal O_X)$ its \emph{irregularity }.

We identify invertible sheaves and Cartier divisors and we use the additive and multiplicative notation interchangeably. Linear equivalence is denoted by $\simeq$. Given two divisors $D_1$ and $D_2$ on $X$, we write $D_1\succeq D_2$ (respectively $D_1\succ D_2$) if the divisor $D_1-D_2$ is (strictly) effective.

\vspace{.5cm}\noindent{\bf Acknowledgments.}  M. Mendes Lopes  was partially supported by  FCT/Portugal  through Centro de An\'alise Matem\'atica, Geometria e Sistemas Din\^amicos (CAMGSD), IST-ID, projects UIDB/04459/2020 and UIDP/04459/2020.  R. Pardini was  partially supported by  project PRIN 2017SSNZAW$\_$004 ``Moduli Theory and Birational Classification"  of Italian MIUR and is a member of GNSAGA of INDAM.   S.Tirabassi was partially supported by the Knut and Alice Wallenberg Foundation under grant no. KAW 2019.0493. 

\section{Preliminaries}
In this section, we provide the background material and preliminary lemmas which are essential for our main result. In particular, we review the logarithmic ramification formula in \ref{LogRam}, and give a quick introduction to quasi-abelian varieties and quasi-Albanese morphism in \ref{sub: qA}. For the sake of brevity, we do not  cover here the general theory of log-varieties and WWPB maps. An exhaustive introduction to the theory of logarithmic forms can be found in \cite{Ii82}. For a quick overview, with also an introduction to WWPB equivalence, we refer to our own \cite{MLPT2022}.
\subsection{Logarithmic ramification formula }\label{LogRam}
Let $V$, $W$  be smooth varieties of dimension $n$ and let  $h\colon V\to W$ be a dominant morphism. Let $g\colon X\to Y$ be a morphism extending $h$, where $X$, $Y$ are smooth compactifications of $V$, respectively $W$, such that  $D:=X\setminus V$ and $\Delta:=Y\setminus W$ are simple normal crossings ({\em snc} for short)  divisors. Then the pull back of a  logarithmic $n$-form on $Y$ is a  logarithmic $n$-form on $X$, and a local computation shows that there is an effective divisor $\ol R_g$ of $X$ - the {\em logarithmic ramification divisor} - such that the following linear equivalence holds:
\begin{equation}\label{eq: logram}
K_X+D\simeq g^*(K_Y+\Delta )+\ol R_g.
\end{equation} 
Equation \eqref{eq: logram} is called the {\em logarithmic ramification formula} (cf. \cite[\S 11.4]{Ii82}). 

 In \cite[Lemma 1.8]{MLPT2022} we observed the following useful fact:
\begin{lm}\label{lem: Rg} In the above set-up, denote by $R_g$ the (usual) ramification divisor of $g$. Let   $\Gamma$ be  an irreducible divisor   such  that $g(\Gamma)\not \subseteq \Delta$.
Then $\Gamma$ is a component of  $\ol R_g$ if and only if $\Gamma$ is a component of $D+R_g$.
\end{lm}

Next we characterize the crepant exceptional divisors, namely those that do not contribute to $\ol R_g$: 
\begin{lm}\label{lm:Lemma2}
 Let $g:X\rightarrow Y$ be a birational morphism of smooth projective varieties, let $D_X\subseteq X$ and $D_Y\subseteq Y$  be snc divisors such that $g(X\backslash D_X)\subseteq Y\backslash D_Y$, and consider $E$ an irreducible $g$-exceptional divisor. If $E$ does not appear in $\overline{R}_g$, then $g(E)$ is an irreducible component of the intersection of some components of $D_Y$. \end{lm}

\begin{proof}
 Let $p\in E$ general and consider $q:=g(p)\in g(E)$.
 
  If $g(E)$ is not contained in $D_Y$, then $q$ is not in $D_Y$ and so, by Lemma \ref{lem: Rg},  $E\preceq \overline{R}_g$. 
  
  Thus, we can assume that $g(E)\subseteq D_Y$ and, therefore, $E\preceq D_X$. Denote by $D_1,\ldots,D_k$ the components of $D_Y$ containing $g(E)$ and let $x_i$ be local equations for $D_i$ near $q$.\par
 Assume by contradiction that $g(E)$ is a proper subset  of $D_1\cap\ldots\cap D_k$. Since $q$ is general,   $g(E)$ is smooth near $q$,  and hence  locally a complete  intersection. Since $D_Y$ is snc,  we can complete $x_1,\ldots,x_k$ to a system of local coordinates $x_1,\ldots,x_k,x_{k+1},\ldots x_n$ such that around $q$ the set  $g(E)$  is locally defined by the vanishing of $x_1,\ldots,x_{k+s}$ for some  $s\ge 1$. In this setting, locally around $q$ the line bundle $\mathcal{O}_Y(K_Y+D_Y)$ is generated by
\begin{align*}
 \sigma&:=\frac{dx_1}{x_1}\wedge\cdots\wedge\frac{dx_k}{x_k}\wedge dx_{k+1}\wedge\cdots\wedge dx_n\\
 &\:=(x_{k+1}\cdots x_{k+s})\frac{dx_1}{x_1}\wedge\cdots\wedge\frac{dx_{k+s}}{x_{k+s}}\wedge dx_{k+s+1}\wedge\cdots\wedge dx_n.
\end{align*}
Choose  a local equation $t$ for $E$ near $p$. Then we can write
$$
g^*x_i=t^{\alpha_i}h_i,
$$
with $h_i$ regular and not divisible by $t$. Observe that $\alpha_i>0$ if $i\leq k+s$, while $\alpha_i=0$ for $i>k+s$. Thus we get 
$$
g^*\sigma =ht^{\alpha_{k+1}+\cdots+\alpha_{k+s}}\cdot\frac{dt}{t}\wedge \tau,
$$
where $h$ is a regular function and $\tau$ is a $(n-1)$-form without poles along $E$. We conclude that $E$ appears in $\overline{R}_g$ with multiplicity at least $\alpha_{k+1}+\cdots+\alpha_{k+s}>0$.  This is a contradiction because we are assuming that $E$ does not appear in $\overline{R}_g$.
\end{proof}

We close the section by proving a partial converse of the previous lemma: 

\begin{lm}\label{lm: newcompactification}
 Let $W$ be a quasi-projective variety with compactification $Y$ and snc boundary $\Delta$. Let $C\subseteq \Delta$ be an irreducible component  of the  intersection of some  components of $\Delta$, and denote by $\epsilon_C:Y_C\rightarrow Y$ the blow-up of $C$.  If $\Delta_C:=\epsilon_C^{-1}\Delta$ is the set-theoretic pre-image of $\Delta$, then $Y_C$ is also a  compactification of $W$ with snc boundary $\Delta_C$ and the logarithmic ramification divisor $\overline{R}_{\epsilon_C}$ is trivial.
\end{lm}
\begin{proof}
 It is clear that $W\simeq Y_C\backslash \Delta_C$ and that $\Delta_C$ has snc support. Thus we need just to verify the statement about the  logarithmic ramification divisor. Let $k$ be the codimension of $C$ in $Y$, and denote by $E$ the exceptional divisor of the blow-up. Then
 \begin{align*}
  K_{Y_C}+\Delta_C\simeq f^*(K_Y)+(k-1)E+f^*(\Delta)-(k-1)E=f^*(K_Y+\Delta).
 \end{align*}
We conclude directly by  the logarithmic ramification formula \eqref{eq: logram} that $\overline{R}_{\epsilon_C}= 0$. 
\end{proof}

\subsection{Quasi-abelian varieties and their compactifications}\label{sub: qA}
A  \emph{quasi-abelian variety} - in some sources also called a \emph{semiabelian variety} - is a connected algebraic group $G$ that is an extension of an abelian variety $A$ by an algebraic torus. More precisely, $G$ sits in the middle of an exact sequence of the form

 \begin{equation}\label{eq:sucQA}1\rightarrow \mathbb{G}_m^r\longrightarrow G\longrightarrow A\rightarrow 0.\end{equation}
 We call $A$ the {\em compact part}  and $\mathbb{G}_m^r$ the {\em linear part} of $G$.
 
We recall the following from \cite[\S 10]{III}:

\begin{prop}\label{prop: compactification}
 Let $G$ be a quasi-abelian variety,  let $A$ be its compact part and let $r:=\dim G-\dim A$. Then there exists a compactification $G\subset Y$ such that:
\begin{itemize}
\item[(a)] $Y$  is a $\mathbb P^r$-bundle over $A$;
\item[(b)] $\Delta:=Y\setminus G$ is a simple normal crossing divisor and  $K_Y+\Delta\simeq 0$;
\item [(c)] the natural $\mathbb{G}_m^r$-action on $G$ extends to $Y$. 

\end{itemize}

\end{prop}
\begin{proof} As explained in \cite[\S 10]{III}, the variety $G$ is a principal $\mathbb{G}_m^r$-bundle over $A$. 
We let   $\rho\colon \mathbb{G}_m^r\to \Aut(\pp^r)$ be the inclusion that maps  $(\lambda_1,\dots  \lambda_r)$ to the automorphism represented by  $\Diag(1,\lambda_1,\dots \lambda_r)$ and take  $Y:= G\times_{\rho}\pp^r$. By construction, $Y=\Proj_A(\OO_A \oplus L_1\oplus\dots  \oplus L_r)$, for some  $L_i\in \Pic(A)$, the boundary $\Delta$ is the sum of the relative hyperplanes $\Delta_0,\Delta_1, \dots \Delta_r$ given by the inclusions of $\OO_Y, L_1, \dots L_r$ into $\OO_A \oplus L_1 \oplus.\oplus L_r$ and the $\mathbb{G}_m^r$-action clearly extends to $Y$. Since $K_A$ is trivial, the formula for the canonical class of a projective bundle reads $K_Y\simeq-(\Delta_0+\dots +\Delta_r)$. \end{proof}
In the sequel, we will denote by $q\colon Y\rightarrow A$ the structure map giving the $\mathbb{P}^r$ bundle structure.
 \subsubsection{The quasi-Albanese map} 
 The classical construction of the Albanese variety of a projective variety can be extended to the non projective case, by replacing regular 1-forms by  logarithmic ones and  abelian varieties by quasi-abelian ones. The key fact is that by Deligne \cite{De71}  logarithmic 1-forms are closed (for the details of the  construction see \cite{III}, \cite[Section 3]{Fu14}).
 As a consequence, for every $V$ smooth quasi projective variety, there exists a quasi-abelian variety $A(V)$ and a morphism $a_V\colon V\to A(V)$ such that any morphism $h\colon V\rightarrow G$  to a quasi-abelian variety  factors through $a_V$ in a unique way up to translation.

\section{Proof of theorem \ref{thm:main}}

Let $V$ be a smooth complex quasi-projective  variety of dimension $n\geq 2$  satisfying   $\overline{\kappa}(V)=0$ and   $\overline{q}(V)=n$.  Let  $A(V)$ be the quasi-Albanese variety of $V$,  and  $a_V\colon V\rightarrow A(V)$ the quasi-Albanese morphism.

Pick  a compactification   $(Y,\Delta)$ of $A(V)$ as in Proposition  \ref{prop: compactification}. Choose then a compactification $(X,D)$ of $V$ such that the quasi-Albanese morphism $a_V\colon V\to A(V)$ extends to a morphism $f\colon X\to Y$.  Recall that $f$ is birational by Kawamata's theorem (\cite[Thm. 28]{Ka81}). 

We observe that proving Theorem A is equivalent  to showing that every irreducible component of $D$ not contracted by $f$ is mapped to  $\Delta$. In fact, let $D'$ be the divisor formed by the irreducible components of $D$ which are not mapped to $\Delta$, assume  that $D'$ is contracted by $f$ and set $Z':=f(D')$.  This is a closed subset of codimension at least 2 in $Y$ and thus $Z:=Z'\cap A(V)$ is a closed subset of codimension at least 2.  Denoting by  $V^0$ the open set $V\backslash a_V^{-1}(Z)$,  we see that the restriction
$$
{a_V}_{|V^0}\colon V^0\rightarrow A(V)\backslash Z
$$
is proper.\par
So we argue by contradiction supposing that  there is an irreducible  component  $H$ of $D$ such that  $f(H)=:\bar H$ is a divisor not contained in $\Delta$.   Then 
\begin{equation}\label{eq:H}
H=f^*(\bar H)-\sum m_i E_i,
\end{equation}
where the divisors $E_i$ are $f$-exceptional divisors. By the logarithmic ramification formula  we have
$$
K_X+D\simeq f^*(K_Y+\Delta)+\overline{R}_f\simeq \overline{R}_f,
$$ 
where the last equality holds because  $K_Y+\Delta \simeq 0$ (see Proposition \ref{prop: compactification}).

We split the proof in two main steps: first, we give an argument in the special case in which all the $E_i$ appearing in \eqref{eq:H} belong to $\overline{R}_f$; afterwards, we show that  the general situation can always be reduced to the special case.\par
\subsection{Special case: all the $E_i$'s are components of $\overline{R}_f$}
Since $H$ is also a component of $\overline{R}_f$ by Lemma \ref{lem: Rg},   we have that 
$$
K_X+D\simeq\overline{R}_f\succeq H+\sum E_i= f^*\overline{H}-\sum(m_i-1)E_i.
$$
Choose $\alpha$ a rational number, $\alpha \in (0,1)$ such that $\alpha\geq \frac{m_i-1}{m_i}$ for every $i$. Then we can write
\begin{align*}
  f^*\overline{H}-\sum(m_i-1)E_i &= \alpha f^*\overline{H}-\sum(m_i-1)E_i+(1-\alpha)f^*\overline{H}\\
  &\succeq \alpha\left( f^*\overline{H}-\sum m_iE_i\right)+(1-\alpha)f^*\overline{H}\\
  &= \alpha H+(1-\alpha)f^*\overline{H}.
\end{align*}
It follows that $\kappa(Y,\overline{H})\le \ol \kappa(V)$. We are going to show that $\kappa(Y,\overline{H})>0$, against the hypothesis $\overline{\kappa}(V)=0$.\par
If $\overline{H}$ is the pullback of an effective divisor $\Gamma$  via the map $q\colon Y\rightarrow A$, then $\kappa(Y,\overline{H})\ge \kappa(A,\Gamma)>0$, because an effective divisor on an abelian variety is the pull back of an ample divisor on a quotient abelian variety. 

Thus we may assume that this is not the case, namely that $\overline{H}$ restricts to a divisor on the general fiber $F$ of the projective bundle $q\colon Y\rightarrow A$. Since $\overline{H}$ is not contained in $\Delta$,  its restriction to the general fiber of $q$ is not contained in $\Delta_{|F}$. In particular $\overline{H}_{|F}$ is not invariant with respect to the $\mathbb{G}_m^r$-action on $F$ (see Proposition  \ref{prop: compactification}). This means that we can find a  $g\in\mathbb{G}_m^r$ such that $g^*\overline{H}_{|F}\neq \overline{H}_{|F}$. Since $g^*\overline{H}$ and $\overline{H}$ are linearly equivalent, we obtain $h^0(Y,\overline{H})\ge 2$ and we are done.
\subsection{Conclusion of the proof}
Up to reordering, we can assume that $E_1\npreceq\overline{R}_f$. By Lemma \ref{lm:Lemma2} we have that $f(E_1)$ is a component of a complete intersection of components of $\Delta$. Let $\epsilon\colon Y_1\rightarrow Y$ be the blow up of $f(E_1)$. Observe that  $Y_1$ is smooth, since $f(E_1)$ is so. By Lemma \ref{lm: newcompactification}, $(Y_1,\Delta_1)$ is a compactification of $A(V)$ with $K_{Y_1}+\Delta_1\simeq 0$. In addition  the $\mathbb{G}_m^r$-action on $A(V)$ extends to $Y_1$ and preserves the components of $\Delta_1$.\par
By the universal property of the blow-up we can factor $f$ as follows:
$$
\xymatrix{X\ar[dr]_{f_1}\ar[rr]^f&&Y\\\
&Y_1\ar[ur]_{\epsilon_1}&
}
$$

Let $\overline{H}_1:=f_1(H)$. This is a divisor not contained in $\Delta_1$, otherwise $f(H)$ would be contained in $\Delta$.  Then
$$f_1^*\overline{H}_1= H+\sum_{i\ge 2}m_i E'_i,$$
where the $E'_i$ are $f_1$-exceptional divisors. If all the $E'_i$ appear  in the  logarithmic ramification divisor $\overline{R}_{f_1}$, we can argue as in the special case and we are done. Otherwise, we can repeat the above construction. As at every stage we provide a blow-up $Y_n$ of $Y_{n-1}$ and a factorization $f=\epsilon_1\circ\cdots \circ \epsilon_n\circ f_n$, after finitely many iterations this process must stop, and again we can argue as in the special case.


\end{document}


\maketitle 

\begin{abstract}
We give an alternative proof of Theorem A in the 
paper:~M.~Mendes Lopes, R.~Pardini, S.~Tirabassi, 
A footnote to a theorem of Kawamata. We 
also explain how to fill a gap in the original proof. 
\end{abstract}

\section{Introduction}

In this paper, we give an alternative proof of the 
following theorem, which is the main result of
\cite{mendes}:

\begin{thm}[{see \cite[Theorem A]{mendes}}]\label{p-thm1.3}
Let $X$ be a smooth variety defined over $\mathbb C$ 
with logarithmic Kodaira 
dimension $\overline \kappa (X)=0$ 
and logarithmic irregularity $\overline q(X)=\dim X$. 
Then the quasi-Albanese map $\alpha \colon X\to A$ 
is birational and there exists a closed subset $Z$ of $A$ 
with $\mathrm{codim}_A Z\geq 2$ such that 
$\alpha\colon X\setminus \alpha^{-1}(Z)\to A\setminus Z$ 
is proper.  
\end{thm}

The proof given in  \cite{mendes} 
contains a gap, noticed by the first 
named author of this paper; in \S \ref{sec: original-proof} we 
explain this gap 
and how to avoid it.  
\begin{ack}
We are grateful to Mark Spivakovsky and 
Beatriz Molina-Samper for informing us 
of the results on embedded resolutions 
of \cite{molina}, that are crucial 
for the argument in \S \ref{sec: original-proof}. \newline
O.~Fujino was partially 
supported by JSPS KAKENHI Grant Numbers 
JP19H01787, JP20H00111, JP21H00974, JP21H04994. M. Mendes Lopes  was partially supported by  FCT/Portugal  through Centro de An\'alise Matem\'atica, Geometria e Sistemas Din\^amicos (CAMGSD), IST-ID, projects UIDB/04459/2020 and UIDP/04459/2020. 
R.~Pardini was partially supported by Project PRIN 
2022BTA242 of italian MUR and is a member of INdAM--GNSAGA. \ S.~Tirabassi was partially supported by project 2023-0387 funded by the \emph{Vetenskapr\aa det}.
\end{ack}

\noindent{\bf Conventions:} We work 
over the field $\mathbb C$ of complex number and  we  use freely  
Iitaka's theory of quasi-Albanese maps and 
logarithmic Kodaira dimension developed in \cite{iitaka1} and \cite{iitaka2} 
(see also \cite{fujino-quasi-alb}). 
 
\section{Proof of Theorem \ref{p-thm1.3}}\label{sec: proof-of-thm}

In this section we  prove Theorem \ref{p-thm1.3}. 
The birationality of $\alpha\colon X\to A$ is a well-known 
theorem by Kawamata 
(see \cite{kawamata-abelian}),  
hence we are going to  prove the existence of the 
desired closed subset $Z$. 
The proof given here uses Kawamata's subadditivity formula in 
\cite{kawamata-master}. 
Before we start the proof of 
Theorem \ref{p-thm1.3}, we note 
the following fact (see also  \cite[Lemma~2.2]{mendes}):

\begin{rem}[Log 
canonical centers]\label{rem: centers}
Let $X$ be a smooth variety and let $\Delta_X$ be a simple 
normal crossings divisor on $X$, so that 
$(X, \Delta_X)$ is log canonical. 
Let $\Delta_X=\sum _{i\in I} \Delta_i$ 
be the irreducible decomposition of $\Delta_X$;
then a closed subset $W$ of $X$  is a log canonical center 
of $(X, \Delta_X)$ if and only if $W$ is 
an irreducible component of $\Delta_{i_1}\cap 
\ldots \cap \Delta_{i_k}$ for some 
$\{i_1, \dots, i_k\}\subset I$. 
When $\Delta_{i_1}\cap 
\ldots \cap \Delta_{i_k}$ is connected for every 
$\{i_1, \dots, i_k\}\subset I$, we say that 
$\Delta_X$ is a {\em{strong normal crossings divisor}} in 
the terminology of \cite{molina}. 
\end{rem}

We now turn to the proof of Theorem \ref{p-thm1.3}. 

\begin{proof}[Proof of Theorem \ref{p-thm1.3}]
Let $\alpha\colon X\to A$ be the quasi-Albanese 
map. By Kawamata's theorem 
(see \cite[Corollary 29]{kawamata-abelian} and 
\cite[Corollary 10.2]{fujino-quasi-alb}), 
we see that $\alpha$ is birational. 
\begin{step}\label{step1}
Let 
\begin{equation*}
\xymatrix{
0\ar[r]& \mathbb G^d_m\ar[r] &A\ar[r]^-\pi&  B\ar[r]& 0 
}
\end{equation*}
be the Chevalley decomposition. 
Then $A$ is a principal $\mathbb G^d_m$-bundle 
over an abelian variety $B$ in the Zariski topology 
(see, for example, 
\cite[Theorems 4.4.1 and 4.4.2]{bcm})
and there is a natural 
completion $\overline \pi\colon \overline A\to B$ of $\pi\colon A\to B$
where $\bar A$ is a 
${\mathbb P}^d$-bundle. We set 
$\Delta_{\overline A}:=\overline A\setminus A$;  
then $\Delta_{\overline A}$ is a simple 
normal crossings divisor on 
$\overline A$, and 
$(\overline A, \Delta_{\overline A})$ is a log canonical pair. 
 
Let $\overline \alpha\colon \overline X\to 
\overline A$ be a compactification of 
$\alpha \colon X\to A$, that is, 
$\overline X$ is a smooth complete 
algebraic variety containing $X$, 
$\Delta_{\overline X}:=
\overline X\setminus X$ is a simple 
normal crossing divisor 
on $\overline X$ and $\overline\alpha$ is a morphism extending $\alpha$.
\begin{claim}
Let $D$ be an irreducible component of 
$\Delta_{\overline X}$ such that 
$\overline \alpha(D)$ is a divisor. 
Then $\overline \pi\colon \overline \alpha (D)\to B$ 
is dominant.  
\end{claim}
\begin{proof}[Proof of Claim] 
We set $D_1:=\overline \alpha (D)$. 
If $\overline \pi\colon D_1\to B$ is not dominant, 
then we can write $D_1={\overline \pi}^*D_2$ 
for some prime divisor $D_2$ on $B$. 
By Remark \ref{rem: centers} every log canonical center 
of $(\overline A, \Delta_{\overline A})$ 
dominates $B$, therefore $D_1$ does not 
contain any log canonical 
center (so in particular it is not a component of $\Delta_{\overline A}$). 
Hence 
\begin{equation*}
K_{\overline X}+\Delta_{\overline X} 
-\overline \alpha^*(K_{\overline A}+\Delta_{\overline A}) 
\geq \varepsilon \overline \alpha^*D_1 
\end{equation*} 
for some $0<\varepsilon \ll 1$ 
since the support of 
$D_1$ does  not contain any  log canonical center of $(\overline A, 
\Delta_{\overline A})$. 
By construction, we have $K_{\overline A}+\Delta_{\overline A}
\sim 0$. 
Thus we obtain  
\begin{equation*}
0=\overline \kappa (X)=\kappa(\overline X, K_{\overline X}
+\Delta_{\overline X})\geq \kappa (\overline X, 
{\overline \alpha}^*D_1)=\kappa (\overline A, D_1)
=\kappa (B, D_2)>0, 
\end{equation*} 
where the last inequality follows from the fact that $D$ is a nonzero effective divisor on the abelian variety $B$.
This  contradiction proves the claim.
\end{proof}
\end{step}

\begin{step}
We assume that there exists an irreducible component 
$D$ of $\Delta_{\overline X}$ such that 
$\overline \alpha (D)$ is a divisor 
with $\overline \alpha(D)\not 
\subset \overline A\setminus A$ and we set 
$D':=\overline \alpha (D)\cap A$. 
By  the  Claim in Step \ref{step1}, $D'$ dominates    $B$,  
therefore 
we can find a subgroup $\mathbb G_m$ of $A$ 
such that $\varphi|_{D'}\colon D'\to A_1$ is dominant, where 
\begin{equation*}
\xymatrix{
0\ar[r]& \mathbb G_m\ar[r] &A\ar[r]^-\varphi&  A_1\ar[r]& 0. 
} 
\end{equation*}  
Note that $A$ is a principal $\mathbb G_m$-bundle over $A_1$ in 
the Zariski topology (see, for example, 
\cite[Theorems 4.4.1 and 4.4.2]{bcm}). 
We take a compactification 
\begin{equation*}
\xymatrix{
f^\dag
\colon X^\dag \ar[r]^-{\alpha^\dag}& A^\dag \ar[r]^-{\varphi^\dag}& A^\dag_1
}
\end{equation*}
of 
\begin{equation*}
\xymatrix{
f
\colon X \ar[r]^-{\alpha}& A\ar[r]^-{\varphi}& A_1, 
}
\end{equation*}
where $X^\dag$, $A^\dag$, and $A^\dag_1$ are smooth complete 
algebraic varieties such that 
$X^\dag\setminus X$, $A^\dag\setminus A$, 
and $A^\dag_1\setminus A_1$ are 
simple normal crossing divisors. 
Let $F$ be a general fiber of $f$. 
Then $\overline \kappa (F)=1$ since $\varphi|_{D'}\colon D'\to A_1$ is 
dominant. 
Note that $A_1$ is a quasi-abelian variety,
hence we have $\overline \kappa (A_1)=0$. 
By Kawamata's theorem (see \cite[Theorem 1]{kawamata-master} and 
\cite[Chapter 8]{iitaka-conjecture}), 
we obtain  
\begin{equation*}
0=\overline \kappa (X)\geq \overline \kappa (F)
+\overline \kappa (A_1)=1. 
\end{equation*} 
This is a contradiction, showing that  
 every irreducible component of $\Delta_{\overline X}$ 
which is not contracted by $\overline \alpha$ is mapped to 
$\Delta_{\overline A}$. 
Let $\Delta'$ be the sum of the irreducible 
components of $\Delta_{\overline X}$ which are not 
mapped to $\Delta_{\overline A}$. 
It is obvious that 
$\mathrm{codim}_{\overline A} \overline \alpha(\Delta')\geq 2$ holds 
since $\Delta'$ is $\overline \alpha$-exceptional. 
We put $Z:=\overline \alpha (\Delta')\cap A$. 
Then $Z$ is a closed subset of $A$ with $\mathrm{codim}_A Z\geq 2$. 
It is easy to see that 
\begin{equation*}
\alpha\colon X\setminus \alpha^{-1}(Z)\to A\setminus Z 
\end{equation*} 
is proper. 
\end{step}
This finishes  the proof of Theorem \ref{p-thm1.3}. 
\end{proof}

\section{On the original proof in \cite{mendes}}\label{sec:  original-proof}

The final part (\S 3.2) of the original proof 
of Theorem \ref{p-thm1.3} given 
in \cite{mendes} contains the 
unsubstantiated  claim that, if 
$\overline X$ contains a divisor 
contracted to a log canonical center 
$W$ of $(\overline A, \Delta_{\overline A})$, 
then $\overline \alpha$ factors though the 
blow-up of $\overline A$ along $W$.  
This claim is then  used to reduce 
the proof to the special case treated in \S 3.1 of \cite{mendes}.

We explain here  how to reduce the 
proof to the special case  by a different argument.
\bigskip

Let $T$ be a smooth variety and let $\Delta_T$ be a simple 
normal crossing divisor on $T$. 
Let $f\colon T'\to T$ be the 
blow-up of a log canonical center  of $(T, \Delta_T)$. 
We set 
\begin{equation*}
K_{T'}+\Delta_{T'}:=f^*(K_T+\Delta_T). 
\end{equation*}  
Then it is easy to see 
that $T'$ is smooth and $\Delta_{T'}$ is a simple normal 
crossing divisor on $T'$. We call $f\colon (T', \Delta_{T'})\to (T, \Delta_T)$ 
a {\em{crepant pull-back}} of $(T, \Delta_T)$. 

Now, with the same notation 
and assumptions of \S \ref{sec: proof-of-thm}, 
let $D$ be a component 
of $\Delta_X$ such that 
$\overline\alpha(D)$ is a divisor not 
contained in $\Delta_{\overline A}$. 
In the terminology of \cite{molina} 
the divisor $\Delta_{\overline A}$ is a 
strong normal crossings divisor (see Remark \ref{rem: centers}). 
Then we are precisely in the 
situation of \S 5.1 of \cite{molina} and by Corollary 2, ibid.,  
there is a  finite sequence of crepant pull-backs 
\begin{equation*}
\xymatrix{
(\overline A, \Delta_{\overline A})=:(T_0, \Delta_{T_0})
&(T_1, \Delta_{T_1})\ar[l]_-{f_1}
& \cdots\ar[l]_-{f_2}
&(T_k, \Delta_{T_k})\ar[l]_-{f_k}
}
\end{equation*}
such that the strict transform of 
$\overline\alpha(D)$ on $T_k$ does not contain 
any log canonical center of $(T_k,\Delta_{T_k})$.
In addition,  the components of $\Delta_{\overline A}$ 
are preserved by the $\mathbb G_m^d$-action 
and therefore it is easy to check that $\mathbb G_m^d$ acts also 
on $T_1$ and the action preserves the components 
of $\Delta_{T_1}$. An inductive argument then shows that 
$\mathbb G_m^d$ acts 
on $T_k$ and  preserves the components of $\Delta_{T_k}$

So, up to replacing $(\overline A, 
\Delta_{\overline A})$ by $(T_k,\Delta_{T_k})$ 
and modifying $\overline X$ accordingly, 
we may assume that $\overline \alpha(D)$ 
does not contain any log canonical center 
and that the $\mathbb G_m^d$-action extends 
to $\overline A$.
We conclude by observing that 
the argument 
in \cite[\S 3.1]{mendes} (the ``special case'') works in this situation. 
